\documentclass[12pt,twoside]{amsart}
\usepackage{amssymb}

\nonstopmode

\numberwithin{equation}{section}
\font\tengothic=eufm10 scaled\magstep 1
\font\sevengothic=eufm7 scaled\magstep 1
\newfam\gothicfam
         \textfont\gothicfam=\tengothic
         \scriptfont\gothicfam=\sevengothic

\newcommand{\Z}{\mathbb{Z}}

\newcommand{\N}{\mathbb{N}}

\newcommand {\PP}{\mathbb{P}}
\newcommand{\xb}{{\bf x}}

\DeclareMathOperator{\codim}{codim}

\DeclareMathOperator{\projdim}{proj\,dim}

\DeclareMathOperator{\rank}{rank}

\DeclareMathOperator{\Soc}{Soc}

\DeclareMathOperator{\dirsum}{\oplus}

\DeclareMathOperator{\pnt}{\raise 0.5mm \hbox{\large\bf.}}

\newcommand{\ffi}{\varphi}
\newcommand{\fm}{{\mathfrak m}}
\newcommand{\im}{\operatorname{im}}
\newcommand{\coker}{\operatorname{coker}}

\newtheorem{theorem}{Theorem}[section]
\newtheorem{lemma}[theorem]{Lemma}
\newtheorem{proposition}[theorem]{Proposition}
\newtheorem{corollary}[theorem]{Corollary}
\newtheorem{conjecture}[theorem]{Conjecture}

\theoremstyle{definition}
\newtheorem{remark}[theorem]{Remark}
\newtheorem{example}[theorem]{Example}

\begin{document}

\title{Extensions of the multiplicity conjecture}

\author{Juan Migliore}
\address{Department of Mathematics, University of Notre Dame,
Notre Dame, IN 46556, USA}
\email{Juan.C.Migliore.1@nd.edu}
\author{Uwe Nagel}
\address
{Department of Mathematics, University of Kentucky,
715 Patterson Office Tower, Lexington, KY 40506-0027, USA}
\email{uwenagel@ms.uky.edu}
\author{Tim R\"omer}
\address{FB Mathematik/Informatik,
Universit\"at Osnabr\"uck, 49069 Osnabr\"uck, Germany}
\email{troemer@mathematik.uni-osnabrueck.de}

\thanks{Part of the work for this  paper was done while the first author was
sponsored by the National Security Agency  under Grant Number MDA904-03-1-0071}

\begin{abstract}
The Multiplicity conjecture of Herzog, Huneke, and Srinivasan states an upper bound for
the multiplicity of any graded $k$-algebra as well as a lower bound for Cohen-Macaulay
algebras. In this note we extend this conjecture in several directions. We discuss when
these bounds are sharp, find a sharp lower bound in case of not necessarily
arithmetically Cohen-Macaulay one-dimensional schemes of 3-space, and we propose an upper
bound for finitely generated graded torsion modules. We establish this bound for torsion
modules whose codimension is at most two.
\end{abstract}

\maketitle
\tableofcontents

\section{Introduction}

Let $R = k[x_1,\dots,x_n]$ be the polynomial ring over the field $k$ with its standard
grading where $\deg x_i = 1$.  Let $N$ be a finitely generated graded $R$-module. We
denote by $e(N)$ the multiplicity of $N$. When $I$ is a saturated ideal defining a closed
subscheme $V \subset {\mathbb P}^{n-1}$, $e(R/I)$ is just the degree, $\deg V$, of $V$.
Consider a minimal free resolution of $N$:
\[
0 \rightarrow \bigoplus_{j=1}^{b_s} R(-d_{s,j}) \rightarrow \dots \rightarrow
\bigoplus_{j=1}^{b_0} R(-d_{0,j})  \rightarrow N \rightarrow 0.
\]
We define $m_i (N)  := \min \{ d_{i,1}, \dots, d_{i,b_i} \}$ and $M_i (N) := \max \{
d_{i,1}, \dots, d_{i,b_i} \}$.  When there is no danger of ambiguity, we simply write
$m_i$ and $M_i$. The module $N$ has a {\em pure resolution} if $m_i = M_i$ for all $i$.

Our first extension of the Multiplicity conjecture of Herzog, Huneke and Srinivasan
includes also a discussion of sharpness:

\begin{conjecture}
\label{conj1} Let $I \subset R$ be a graded ideal,  $c=\codim R/I$ and assume that $R/I$
is Cohen-Macaulay. Then
$$
\frac{1}{c!} \prod_{i=1}^c m_i \leq e(R/I) \leq \frac{1}{c!} \prod_{i=1}^c M_i
$$
with equality below (resp.\ above) if and only if $R/I$ has a pure resolution.
\end{conjecture}

Indeed, the original conjecture consisted only of the two inequalities, but the statement
about equality is motivated by the results of \cite{MNR}, where this extended conjecture
is shown for Cohen-Macaulay algebras of codimension two and Gorenstein algebras of
codimension 3.

There has been a tremendous amount of activity towards establishing important special
cases of the Multiplicity conjecture (cf., e.g., \cite{HM}, \cite{HS}, \cite{srinivasan}
 \cite{gold}, \cite{GSS}, \cite{roemer}, \cite{GV},  \cite{MNR}, \cite{F} \cite{rosapaper},
\cite{H}, and the survey paper \cite{FS}). In this paper, our goal is to give several
extensions of this conjecture and to prove some initial cases.  First, we discuss
non-Cohen-Macaulay  $k$-algebras.

It is conjectured in \cite{HS} that the upper bound in Conjecture \ref{conj1} is also
true when $R/I$ is not Cohen-Macaulay. However, it is known that the inequality
$\frac{1}{c!} \prod_{i=1}^c m_i \leq e(R/I)$ does not hold in general.  A simple example
is given by the coordinate ring $B$ of two skew lines in ${\mathbb P}^3$, where $e(B) =
2$ but $m_1(B) = 2, m_2(B) = 3$. Thus, our first goal in this note is to examine what it
is that prevents the lower bound from holding, at least in the case of one-dimensional
curves in ${\mathbb P}^3$.  In the process, we find a lower bound. Furthermore, we get an
improved upper bound.

Specifically, let $M(C)$ denote the Hartshorne-Rao module of $C$ (see the beginning of
Section \ref{lower bd P3} for the  definition).  $M(C)$ has been studied a great deal in
the literature, most successfully in work of Rao for liaison of curves in ${\mathbb P}^3$
(see e.g.\ \cite{rao}).  However, in our case is it a {\em submodule} of $M(C)$ that is
important.  Namely, if $A$ is the ideal generated by two general linear forms, then we
denote by $K_A$ the submodule of $M(C)$ annihilated by $A$.  This has finite length
regardless of whether or not $M(C)$ does.  It has been studied in \cite{london} and in
\cite{migbook}, but this is a new application.  We show in Theorem \ref{thm-bounds-curve}
that
\[
\frac{1}{4} m_1(C) m_2(C) - \dim_k K_A \leq \deg C \leq \frac{1}{2} M_1(C)M_2(C) - \dim_k
K_A.
\]

Then we show in Section \ref{sec-bdl} how basic double G-links can be used to approach
the Multiplicity conjecture, also in the non-Cohen-Macaulay case. As applications, we use
this rather general method to provide a unified framework for giving simple proofs that
the Multiplicity conjecture is preserved under regular hypersurface sections and that the
Multiplicity conjecture is true  for standard determinantal ideals. The latter result has
been independently obtained, first by Mir\'o-Roig \cite{rosapaper} and then by Herzog and
Zheng \cite{H}. The hypersurface section result is also contained in \cite{H}.

 As further
application of the methods in Section \ref{sec-bdl} we derive a closed formula for the
degree of any standard determinantal ideal (Theorem \ref{thm-degree-form}).

Note that in \cite{H} Herzog and Zheng have shown that virtually in all cases where the
multiplicity  bounds were previously known, the bounds are sharp if and only if the
algebra $R/I$ is Cohen-Macaulay and has a pure resolution.

Finally, we discuss the case of modules. We propose the following extension of the
Multiplicity conjecture for cyclic modules in \cite{HS}:

\begin{conjecture}
\label{conj2} Let $N$ be a finitely generated graded torsion $R$-module of codimension
$c= n - \dim N$. Then we have that
$$
e(N) \leq \frac{1}{c!} \prod_{i=1}^c (M_i - m_0)
$$
with equality if and only if $N$ is Cohen-Macaulay and has a pure resolution.
\end{conjecture}

Extending a result by Herzog and Srinivasan in \cite{HS}, we  first show this conjecture
for every Cohen-Macaulay module of any codimension $c > 0$ with a quasi-pure resolution,
i.e.\ $m_i \geq M_{i-1}$ for all $i \geq 1$.  This is carried out in Section
\ref{sec-quasi-pure}. The result (cf.\ Theorem \ref{modulequasi}) is used as one
ingredient in Section \ref{sec-mod-2} where we prove one of the main results of this
paper: Conjecture \ref{conj2} is true for every module whose codimension is at most two
(cf.\ Theorem \ref{modulecodim2}). This extends the result for such cyclic modules in
\cite{roemer}. It also provides our strongest evidence for Conjecture \ref{conj2}.

We hope that our conjectured extensions of the Multiplicity conjecture and the methods of
this paper will stimulate further investigations.


\section{A lower bound for non-arithmetically Cohen-Macaulay curves in ${\mathbb
P}^3$} \label{lower bd P3}

Let $C \subset {\mathbb P}^3$ be a one-dimensional scheme; we do not necessarily
assume that the saturated ideal  $I_C \subset R = k[x_1,x_2,x_3,x_4]$ is unmixed.  Let
\[
M(C) := \bigoplus_{t \in {\mathbb Z}} H^1({\mathbb P}^3, {\mathcal I}_C(t))
\]
 be the {\em deficiency module} of $C$ (also known as the {\em Hartshorne-Rao
module} of $C$).  $C$ is arithmetically Cohen-Macaulay if and only if $M(C) = 0$, and $C$
is locally Cohen-Macaulay and equidimensional (i.e.\ $I_C$ is unmixed) if and only if
$M(C)$ has finite length.  If $L$ is a general linear form, then $L$ induces a
multiplication $(\times L) : M(C) \rightarrow M(C)(1)$.  The kernel, $K_L = [0:_{M(C)}
(L)]$, has finite length regardless of whether or not $M(C)$ does. To simplify notation
we write for any ideal $I \subset R$ throughout the paper, $m_i (I)$ and $M_i (I)$ for
$m_i (R/I)$ and $M_i (R/I)$ and similarly $m_i (C)$ and $M_i (C)$ for $m_i (I_C) = m_i
(R/I_C)$ and $M_i (I_C) = M_i (R/I_C)$, respectively. Since in this section we will
compare resolutions of modules and ideals, we will use for an $R$-module $N$ (in contrast
to the notation in the introduction) $m_i (N)$ and $M_i (N)$ to refer to the $(i-1)$-st
syzygy module of $N$.

If $A$ is the ideal generated by two general linear forms, then we set $K_A$ to
be the submodule of $M(C)$ annihilated by $A$.  In \cite{london}, it was shown
that $K_A$ plays a very interesting role in the study of one-dimensional
subschemes of ${\mathbb P}^3$.  For instance, a famous result of Dubreil for
codimension two arithmetically Cohen-Macaulay curves in ${\mathbb P}^3$ (in fact
for codimension two arithmetically Cohen-Macaulay subschemes of any projective
space) says that
\[
\nu(C) \leq m_1(C) +1,
\]
where $\nu(C)$ is the number of minimal generators of the homogenous ideal $I_C$.  This
result is not true for non arithmetically Cohen-Macaulay curves.  However, using $K_A$
allows one to generalize Dubreil's theorem:
\[
\nu(C) \leq m_1(C) + 1 + \nu(K_A).
\]
Here $C$ can in fact be {\em any} subscheme of projective space of dimension $\leq 1$.

The main purpose of this section is to show how $K_A$ can be used to similarly extend the
lower bound of the Multiplicity conjecture for non-arithmetically Cohen-Macaulay
one-dimensional schemes.  In this geometric context we prefer to write $\deg C$ for
$e(R/I_C)$.  Recall that the standard formulation is known to be false: for instance, if
$C$ is the union of two skew lines then $\frac{1}{2} m_1(C) m_2(C) = \frac{1}{2}(2)(3) =
3 > 2 = \deg C$.  In \cite{MNR} it was shown that even a bound of the form $\frac{1}{p}
m_1(C) m_2(C) \leq \deg C$ is false, for any $p \geq 1$ (see Example \ref{sharp ex from
MNR} below), even if we were to restrict to unmixed one-dimensional schemes (which in
fact we will not have to do).  Here we show that $K_A$ is the right ``correction term."
The main result of this section is:

\begin{theorem} \label{thm-bounds-curve}
Let $C \subset {\mathbb P}^3$ be a one-dimensional scheme.Then  we have
\[
\frac{1}{4} m_1(C) m_2(C) - \dim_k K_A \leq \deg C \leq \frac{1}{2} M_1(C) M_2(C) -
\dim_k K_A.
\]
\end{theorem}

The proof is based on the following slightly technical result.

\begin{lemma} \label{lowerbd-curves}
Let $C \subset {\mathbb P}^3$ be a one-dimensional scheme with saturated ideal
$I_C$.  Let
$A = (L_1, L_2)$ be the ideal generated by two general linear forms in $R =
k[x_0,x_1,x_2,x_3]$.  Let $J = \frac{I_C + A}{A} \subset R/A := T  \cong
k[x,y]$. Let
\[
m_2 = \min \{ m_2(C), m_1(K_A) +2 \}.
\]
Then
\[
\frac{1}{2} m_1(C) m_2 - \dim_k K_A \leq \deg C \leq \frac{1}{2} M_1(C)M_2(C) - \dim_k
K_A.
\]
\end{lemma}

\begin{proof}

From \cite{migbook} page 49 we know that
\[
K_A(-2) \cong \frac{I_C \cap A}{I_C \cdot A}.
\]
Now consider the exact sequence
\begin{equation} \label{ses}
\begin{array}{cccccccccccccccc}
0 & \rightarrow & \displaystyle \frac{I_C \cap A}{I_C \cdot A} & \rightarrow &
\displaystyle \frac{I_C}{I_C
\cdot A} & \rightarrow & \displaystyle \frac{I_C}{I_C \cap A} &
\rightarrow & 0.  \\
&& || &&&& || \\
&& K_A(-2) && && \displaystyle \frac{I_C + A}{A}
\end{array}
\end{equation}
The latter module, $\frac{I_C + A}{A} =: J$, is isomorphic to an ideal in
$k[x_0,x_1,x_2,x_3]/A \cong k[x,y] =: T$.  As such, it satisfies the Multiplicity
conjecture.  Observe that the graded Betti numbers of the $T$-module $\frac{I_C}{I_C
\cdot A}$, over $T$, are the same as those of $I_C$ over $R$.

Considering resolutions over $T$, we get the diagram
\begin{equation} \label{diag}
\begin{array}{cccccccccccccccc}
&& 0 && 0 \\
&& \downarrow && \downarrow \\
&& {\mathbb K}_3 (-2) && {\mathbb F}_3 && 0 \\
&& \downarrow && \downarrow && \downarrow \\
&& {\mathbb K}_2 (-2) &&  {\mathbb F}_2 && {\mathbb G}_2 \\
&& \downarrow && \downarrow && \downarrow \\
&& {\mathbb K}_1(-2) && {\mathbb F}_1 && {\mathbb G}_1 \\
&& \downarrow && \downarrow && \downarrow \\
0 & \rightarrow & K_A(-2) & \rightarrow & \displaystyle \left ( \frac{I_C}{A
\cdot I_C}
\right ) & \rightarrow & J & \rightarrow & 0 \\
&& \downarrow && \downarrow && \downarrow \\
&& 0 && 0 && 0
\end{array}
\end{equation}
The mapping cone immediately gives that
\[
{\mathbb K}_3 (-2) \cong {\mathbb F}_3, \hbox{ and ${\mathbb K}_3$ splits in the
mapping cone}
\]
and that
\begin{equation} \label{summand}
{\mathbb K}_2 (-2) \hbox{ is a direct summand of ${\mathbb F}_2$ and splits in
the mapping cone}.
\end{equation}

\bigskip

Now, using (\ref{ses}), we have for any $t$
\[
-\dim_k  J_t = - \dim_k \left ( \frac{I_C}{A \cdot I_C} \right )_{t} + \dim_k (K_A)_{t-2}
\]
so
\[
\begin{array}{rcl}
\dim_k (T/J)_t & = & \Delta^2 \dim_k (R/I_C)_t  + \dim_k (K_A)_{t-2} \\ \\
& = & [\Delta \dim_k (R/I_C)_t] - [\Delta \dim_k (R/I_C)_{t-1}] + \dim_k (K_A)_{t-2}.
\end{array}
\]
Now, we know that for $t \gg 0$ we have $\Delta \dim_k (R/I_C)_t = \deg C$ and $\dim_k
(K_A)_t = 0$, while for $t < 0$ we have $\Delta \dim_k (R/I_C)_t = \dim_k (K_A)_t = 0$.
Now let $t_0 \gg 0$ and sum over all $t \leq t_0$:
\begin{equation} \label{important computation}
\begin{array}{rcl}
\deg C & = & \Delta \dim_k (R/I_C)_{t_0} \\
& = & \displaystyle \sum_{t=0}^{t_0} [\Delta \dim_k(R/I_C)_t - \Delta \dim_k
(R/I_C)_{t-1}]
\\ & = & \displaystyle \sum_{t=0}^{t_0} \dim_k (T/ J)_t - \dim_k (K_A)_{t-2} \\
& = & e (T/ J) - \dim_k K_A .
\end{array}
\end{equation}
Now, as noted above, $ J \subset T \cong k[x,y]$ satisfies the Multiplicity conjecture.
Note also that from the definition of $J$ it follows immediately that
\[
m_1(J) = m_1\left ( \frac{I_C}{A \cdot I_C} \right ) = m_1(C),
\]
and from the mapping cone associated to (\ref{diag}) we have
\[
m_2(J) \geq m_2.
\]
  Now we get
\[
\begin{array}{rcl}
\deg C & =  & e (T/ J) - \dim_k K_A \\[2pt]
& \geq & \frac{1}{2} m_1( J) m_2( J) - \dim_k K_A \\[2pt]
& \geq & \frac{1}{2} m_1(C) m_2 - \dim_k K_A
\end{array}
\]
proving the lower bound.

For the upper bound, the same argument works as above (reversing the inequalities and replacing $m_i$ by $M_i$), once we have shown that
\[
M_2(J) \leq M_2(C).
\]
To see this, the mapping cone construction for (\ref{diag}) shows that
\[
M_2(J) \leq \max \{ M_2(C), M_1(K_A(-2)) \}.
\]
But using (\ref{summand}), we see that
\[
M_2(C) \geq M_2(K_A(-2)) > M_1(K_A(-2)),
\]
so we are done.
\end{proof}

\begin{corollary} \label{cor-curves}
With the notation of Lemma \ref{lowerbd-curves}, we furthermore let $H$ be the plane in
${\mathbb P}^3$ defined by $L_1$ and let $C \cap H$ be the geometric hyperplane section
(i.e.\ it is defined by the saturation of $\frac{I_C + (L_1)}{(L_1)}$ in $R/(L_1)$, which
we denote by $I_{C \cap H}$). Then:
\begin{itemize}
\item[(a)] If $m_2(C) \leq m_1(K_A)+2$ then
\[
\frac{1}{2} m_1(C) m_2(C) - \dim_k K_A \leq \deg C .
\]

\item[(b)] If $m_1(K_A) +2 < m_2(C)$ then the following both hold:

\begin{itemize}

\item[(i)] $\displaystyle \frac{1}{2} m_1(C) m_2(C \cap H) - \dim_k K_A \leq \deg C .$ \\

\item[(ii)] $\displaystyle \frac{1}{4} m_1(C) m_2(C) - \dim_k K_A \leq \deg C $.
\end{itemize}
\end{itemize}
\end{corollary}

\begin{proof}
Part (a) is immediate from Theorem \ref{lowerbd-curves}.  For (b), we first
prove (i). Consider the exact sequence
\[
0 \rightarrow K \rightarrow K_{L_1}(-1) \stackrel{\times
L_2}{\longrightarrow} K_{L_1} \rightarrow B \rightarrow 0,
\]
where $K$ is the kernel and $B$ is the cokernel, respectively, of the multiplication by
$L_2$. Note that $K = K_A(-1)$.  The Socle Lemma (\cite{HU}, Corollary 3.11) then gives
that $m_1(K)
> m_1(\Soc(B))$.  Now, consider the commutative diagram
\[
\begin{array}{cccccccccccccccccccc}
&&&&&&&&0 \\
&&&&&&&& \downarrow \\
&& 0 && 0 && 0 && K_A(-2) \\
&& \downarrow && \downarrow && \downarrow && \downarrow \\
0 & \rightarrow & I_C(-2) & \stackrel{\times L_1}{\longrightarrow} & I_C(-1) &
\rightarrow & I_{C \cap H}(-1) & \rightarrow & K_{L_1}(-2) & \rightarrow & 0\\
&& \phantom {\times L_2} \downarrow { \times L_2} && \phantom {\times
L_2} \downarrow { \times L_2} && \phantom {\times L_2} \downarrow {
\times L_2} && \phantom {\times L_2} \downarrow { \times L_2} \\
0 & \rightarrow & I_C(-1) & \stackrel{\times L_1}{\longrightarrow} & I_C &
\rightarrow & I_{C \cap H} & \rightarrow & K_{L_1} (-1) & \rightarrow & 0 \\
&& \downarrow && \downarrow && \downarrow && \downarrow \\
0 & \rightarrow & J(-1) & \stackrel{\times L_1}{\longrightarrow} & J &
\rightarrow & \displaystyle \frac{I_{C \cap H}}{L_2 \cdot I_{C\cap H}} &
\rightarrow & B & \rightarrow & 0 \\
&& \downarrow && \downarrow && \downarrow && \downarrow \\
&& 0 && 0 && 0 && 0
\end{array}
\]
where $J = I_C / (L_2 \cdot I_C)$ and $\frac{I_{C \cap H}}{L_2 \cdot I_{C\cap
H}}$ is isomorphic to an ideal, $J'$, in $T := R/(L_1,L_2)$ with the same graded
Betti numbers as $I_{C \cap H}$ in $R/(L_1)$.  We denote by $\Soc (C \cap H)$
the socle of $T/J'$.   We get that
\[
\begin{array}{rcl}
m_2 & = & m_1(K_A) +2 \\
& = & m_1(K) +1 \\
& \geq & m_1(\Soc(B)) +2 \\
& \geq & m_1(\Soc(C \cap H)) +2 \\
& = & m_2(C \cap H).
\end{array}
\]
Then Theorem \ref{lowerbd-curves} gives the desired result.

For (ii), recall that $K_A(-2) \cong \frac{I_C \cap A}{I_C \cdot A}$.  Let $d =
m_1(K_A(-2)) \geq m_1(C)$.  Let $F \in I_C$ be a generator of minimal degree,
$m_1(C)$.  We claim that there is an element $G \in (I_C)_d$ that is not a
multiple of $F$.  Indeed, since $A$ corresponds to a general line, we can assume
that $F$ does not vanish identically on this line.  If every element of
$(I_C)_d$ were a multiple of $F$ then in particular any element of $(I_C \cap
A)_d$ is of the form $HF$, where $H \in A$, contradicting the fact that $K_A(-2)
= \frac{I_C \cap A}{I_C \cdot A}$ begins in degree $d$.  It follows, since $F$
and $G$ have at worst a Koszul syzygy and $G$ may not itself be a minimal
generator, that
\[
m_2(C) \leq m_1(C) + m_1(K_A(-2)) \leq 2 m_1(K_A(-2)).
\]
Now, in Theorem \ref{lowerbd-curves} we set $m_2 = \min\{m_2(C), m_1(K_A) +2
\}$, which in our current hypothesis is equal to $m_1(K_A(-2))$.  We showed
that
\[
\deg C \geq \frac{1}{2} m_1 (C) \cdot m_2 - \dim_k K_A.
\]
Now we obtain
\[
\begin{array}{rcl}
\deg C & \geq & \displaystyle \frac{1}{2} m_1(C) \left ( \frac{1}{2}
m_2(C)  \right ) - \dim_k K_A \\
& = & \displaystyle \frac{1}{4} m_1(C) m_2(C) - \dim_k K_A
\end{array}
\]
as claimed.
\end{proof}

Summing up we get:

\begin{proof}[Proof of Theorem \ref{thm-bounds-curve}]
This is an immediate consequence of Lemma \ref{lowerbd-curves} and Corollary
\ref{cor-curves}.
\end{proof}

A very strong generalization of the following result can be found in Theorem
\ref{modulecodim2}, but  we include it here as an interesting special case.

\begin{corollary}
Let $C \subset {\mathbb P}^3$ be a one-dimensional scheme.  Then
\[
\deg C \leq \frac{1}{2} M_1(C)M_2(C),
\]
with equality if and only if $C$ is arithmetically Cohen-Macaulay with pure resolution.
\end{corollary}

\begin{proof}
It follows immediately from \cite{MNR} Corollary 1.3, and the observation that $K_A = 0$ if and only if $C$ is arithmetically Cohen-Macaulay.
\end{proof}

\begin{remark}
The analogous statements for the lower bound are not true.  First, it is well known that
\[
\frac{1}{2} m_1(C) m_2(C) \leq \deg C
\]
is false in general (consider for instance two skew lines).  Furthermore, if
\[
\frac{1}{2} m_1(C) m_2(C) = \deg C,
\]
it does not follow that $C$ is arithmetically Cohen-Macaulay.  A simple example is a
rational quartic curve in ${\mathbb P}^3$.
\end{remark}

\begin{corollary}
Let $C \subset {\mathbb P}^3$ be a curve, i.e.\ $I_C$ is an unmixed
one-dimensional scheme (so $\dim_k M(C) < \infty$).  Then
\[
\frac{1}{4} m_1(C) m_2(C) - \dim_k M(C) \leq \deg C.
\]
\end{corollary}

Note, though, that in general $\dim_k M(C)$ can be much greater than $\dim_k K_A$. See
the following example.

\begin{example}
In Theorem \ref{lowerbd-curves}, if it were true that
\[
m_2( J) \geq m_2 \left (
\frac{I_C}{A \cdot I_C} \right ) = m_2(C)
\]
then in fact we would be able to prove that the first assertion of Corollary
\ref{cor-curves} always holds:
\[
\frac{1}{2} m_1(C) m_2(C) \leq \deg C + \dim_k K_A.
\]
Unfortunately, we now show that these assertions are both false in general.

Let $C$ be the disjoint union of two complete intersections of type $(16, 16)$.  Then
$\dim_k M(C) = 65,536$, $\dim_k K_L = 2,736$ (where $K_L$ is the submodule of $M(C)$
annihilated by one general linear form) and $\dim_k K_A = 171$.  We also have $\deg C =
512$, and one can compute
\[
\begin{array}{rclcrclcrclcrcl}
m_1(C) & = & 32 && m_1(K_A) +2 & = & 42 && m_1(J) & = & 32 && m_1(C \cap H) =
31  \\
m_2(C) & = & 48 & & m_2(K_A) +2 & = & 48 && m_2(J) & = & 42 && m_2(C \cap H) =
33 \\
m_3(C) & = & 64 & & m_3(K_A) +2 & = & 64. \\
\end{array}
\]
In particular, $m_2(J) < m_2(C)$ and
\[
768 = \frac{1}{2} m_1(C) m_2(C) > \deg C + \dim_k K_A = 683.
\]
\end{example}

\begin{example} \label{sharp ex from MNR}
In \cite{MNR}, Remark 2.4, we began the discussion of what might be the right
approach to a lower bound for the multiplicity for non-arithmetically
Cohen-Macaulay space curves.  We noted that a bound of the form
\[
\frac{1}{p} m_1(C) m_2(C) \leq \deg C
\]
is not possible, for any fixed value of $p \geq 1$.  We gave as example the curve $C$
with saturated ideal
\begin{equation} \label{eqn of C}
I_C = (x_0,x_1)^t + (F),
\end{equation}
where $F$ is smooth along the line defined by $(x_0,x_1)$ and $d := \deg F \geq t+1$, and
we noted that $\deg C = t$,  $m_1(C) = t$ and $m_2(C) = t+1$.  Hence we would need $p
\geq t+1$, and obviously this can be made arbitrarily large.

Since this example was used to illustrate the difficulty of finding a nice lower bound,
we would like to remark here that our Corollary \ref{cor-curves} fits nicely with this
example.  In Proposition \ref{sharp ex} we will show that a curve of this form {\em
always} gives equality in Corollary \ref{cor-curves} (a).  First, though, we give a
specific numerical example to illustrate not only the sharpness, but the huge difference
that is possible between $\dim_k M(C)$ and $\dim_k K_A$ even in the context of Corollary
\ref{cor-curves} (a).

Specifically, let $t=12$ and $\deg F = 15$. Using {\tt macaulay} \cite{macaulay}, we have
computed $\dim_k M(C) = 56,056$, but $\dim_k K_A = 66$.  Furthermore, $m_1(K_A) = 13,
m_2(K_A) = 24$ and $m_3(K_A) = 25$ (over $T$).  Then the hypothesis of (a) holds, and in
fact we have $\deg C = 12$ and
\[
78 = \frac{1}{2} m_1(C) m_2(C) = \deg C + \dim_k K_A.
\]
\end{example}

\begin{proposition} \label{sharp ex}
Let $C$ be the non-reduced curve with ideal $I_C$ given in (\ref{eqn of C}).
Then
\[
\frac{1}{2} m_1(C) m_2(C) = \deg C + \dim_k K_A.
\]
\end{proposition}

\begin{proof}
We have from (\ref{important computation}) that
\[
\deg C = e( T/J) - \dim_k K_A,
\]
where $A = \frac{I_C + A}{A}$.  In our case, without loss of generality we may
choose $A = (x_2,x_3)$.  Note that then
\[
\begin{array}{rcl}
J & = & \displaystyle \frac{(x_0,x_1)^t + (x_2,x_3) + (F)}{(x_2,x_3)}  \\
& = & (x_0,x_1)^t
\end{array}
\]
viewed in the ring $T = k[x_0,x_1,x_2,x_3]/(x_2,x_3)$.  (Note that in \cite{MNR},
Remark 2.4, we could even have taken $\deg F = t$ and we would still have this
equality.)  Hence we know that $e(T/ J) =
\binom{t+1}{2}$, so we compute
\[
\dim_k K_A = e(T/ J) - \deg C = \binom{t+1}{2} - t = \binom{t}{2}.
\]
On the other hand, we know that $m_1(C) = t, m_2(C) = t+1$.  Indeed, the entire
minimal free resolution of $I_C$ can be computed from a mapping cone using the
exact sequence
\[
0 \rightarrow (x_0,x_1)^{t-1}(-d) \rightarrow (x_0,x_1)^t \oplus R(-d)
\rightarrow I_C \rightarrow 0.
\]
We then immediately see that
\[
\begin{array}{rcl}
\displaystyle \frac{1}{2} m_1(C) m_2(C) & = & \displaystyle \binom{t+1}{2}
\\ & = & \displaystyle t + \binom{t}{2} \\
& = & \deg C + \dim_k K_A
\end{array}
\]
as desired.
\end{proof}


\section{Basic double G-links} \label{sec-bdl}

In this section we prove that under suitable assumptions in arbitrary codimension,  basic
double G-links preserve the property of satisfying the upper and lower bounds of the
Multiplicity conjecture. We use this to provide a unified framework for some results that
also have been independently obtained by Mir\'o-Roig \cite{rosapaper} and by Herzog and
Zheng \cite{H}. Furthermore, we establish a formula for the degree of any standard
determinantal ideal. We believe that this framework will provide further applications.

 We first recall the set-up.

 \begin{proposition}[\cite{KMMNP} Lemma 4.8, Proposition 5.10] \label{KMMNP result}  Let
 $I \subset J$ be homogenous ideals of $R = k[x_0,\dots,x_n]$ such that
 $\hbox{codim } I + 1 = \hbox{codim } J = c+1$.  Let $L \in R$ be a form of degree $d$ such that
 $I : L = I$.  Let $J_1 = I + L \cdot J$.   Then we have

 \begin{itemize}
 \item[(i)] $e(R/ J_1) = d \cdot e(R/ I) + e(R/ J)$.

 \item[(ii)] If $R/I$ satisfies property $G_0$ (Gorenstein in codimension 0) and if $J$ is
 unmixed, then $J$ and $J_1$ are Gorenstein linked in two steps.

 \item[(iii)] We have a short exact sequence
 \[
 0 \rightarrow I(-d) \rightarrow J(-d) \oplus I \rightarrow J_1 \rightarrow 0,
 \]
 where the first map is given by $F \mapsto (LF,F)$ and the second map  is given by
 $(F,G) \mapsto F-LG$.

 \end{itemize}
 \end{proposition}

 \begin{remark}
 Statement (ii) of Proposition \ref{KMMNP result} explains why this process is called Basic
 Double G-linkage.
 However, we will not need this fact here.
 \end{remark}

Recall our convention from Section \ref{lower bd P3} that we write for an ideal $I
\subset R$, $m_i (I)$ for $m_i (R/I)$ and similarly $M_i (I)$ for $M_i (R/I)$.

\begin{theorem} \label{upper bd for BDL}
Under the assumptions of Proposition \ref{KMMNP result}, assume that $I$ and $J$  both
 satisfy the upper bound of the Multiplicity conjecture.
Assume further that
\begin{equation} \label{upper bd assumption}
M_i(J_1) \geq M_i(J)+d, \ \ \ 1 \leq i \leq c+1
\end{equation}
and
\begin{equation} \label{fact2}
M_i(J_1) \geq M_{i-1}(I)+d, \ \ \ 2 \leq i \leq c+1.
\end{equation}
 Then $J_1$ also satisfies the upper bound of the Multiplicity conjecture.
\end{theorem}

\newcounter{tempp}
\setcounter{tempp}{\value{equation}}

\begin{proof}
We have the exact
sequence
\begin{equation} \label{short exact}
0 \rightarrow I(-d) \rightarrow J(-d) \oplus I \rightarrow J_1
\rightarrow 0.
\end{equation}
Let $\mathbb F_\bullet$ be a minimal free resolution for $I$ and
$\mathbb G_\bullet$ be a minimal free resolution for $J$. Consider the exact diagram
\begin{equation} \label{big diagram}
\begin{array}{ccccrclccccccccccccccccccccccc}
&&&&&0 \\
&&&&& \downarrow \\
&& 0 && {\mathbb G}_{c+1}(-d) \\
&& \downarrow &&& \downarrow \\
&& {\mathbb F}_c(-d) && {\mathbb G}_c(-d) & \oplus & {\mathbb F}_c \\
&& \downarrow &&& \downarrow \\
&& {\mathbb F}_{c-1}(-d) && {\mathbb G}_{c-1}(-d) & \oplus & {\mathbb F}_{c-1} \\
&& \downarrow &&& \downarrow \\
&& \vdots &&& \vdots \\
&& \downarrow &&& \downarrow \\
&&{\mathbb F}_3(-d) && {\mathbb G}_3(-d) & \oplus & {\mathbb F}_3 \\
&& \downarrow &&& \downarrow \\
&& {\mathbb F}_2(-d) && {\mathbb G}_2(-d) & \oplus & {\mathbb F}_2 \\
&& \downarrow &&& \downarrow \\
&& {\mathbb F}_1(-d) && {\mathbb G}_1(-d) & \oplus & {\mathbb F}_1 \\
&& \downarrow &&& \downarrow \\
0 & \rightarrow & I(-d) & \rightarrow & J(-d) & \oplus & I & \rightarrow & J_1 & \rightarrow & 0 \\
&& \downarrow &&& \downarrow \\
&& 0 &&& 0

\end{array}
\end{equation}

From this diagram with the induced maps, the  mapping cone gives us the following free resolution
for
$J_1$:
\begin{equation} \label{resolution}
0 \rightarrow
\begin{array}{c}
\mathbb F_c(-d) \\
\oplus \\
\mathbb G_{c+1}(-d)
\end{array}
\rightarrow
\begin{array}{c}
\mathbb F_{c-1}(-d) \\
\oplus \\
\mathbb G_{c}(-d) \\
\oplus \\
\mathbb F_{c}
\end{array}
\rightarrow
\begin{array}{c}
\mathbb F_{c-2}(-d) \\
\oplus \\
\mathbb G_{c-1}(-d) \\
\oplus \\
\mathbb F_{c-1}
\end{array}
\rightarrow
\dots
\rightarrow
\begin{array}{c}
\mathbb F_{1}(-d) \\
\oplus \\
\mathbb G_{2}(-d) \\
\oplus \\
\mathbb F_{2}
\end{array}
\rightarrow
\begin{array}{c}
\mathbb F_1 \\
\oplus \\
\mathbb G_1(-d)
\end{array}
\rightarrow J_1 \rightarrow 0
\end{equation}
For our
purposes we only need information about the largest shift in any free module.

We first note that our hypotheses (\ref{upper bd assumption}) and (\ref{fact2}) are
implied by the following statement:

\medskip

\noindent (\ref{upper bd assumption}$'$) \ \ \ \parbox{5.5 in} {\em In the $i$-th free
module in the  resolution (\ref{resolution}), there is at least one summand of ${\mathbb
G}_i(-d)$ of (total) degree $M_i(J)+d$ that is not split off.}

\medskip

\noindent
This is a very natural situation, but it is not {\em always} true, as seen in Example
\ref{all top deg split} below.

An examination of the diagram (\ref{big diagram})  and the corresponding mapping cone
(including the maps) reveals that the only possible splitting in (\ref{resolution}) comes
in canceling a term $R(-t-d)$ of $\mathbb F_{i}(-d)$ in the $(i+1)$-st free module with a
corresponding one of $\mathbb G_i(-d)$ in the $i$-th free module.  This shows right away
that
\begin{equation}\label{fact1}
M_i(J_1) \geq M_i(I), \ \ \ 1 \leq i \leq c,
\end{equation}
since no summand of $\mathbb F_i$ in the $i$-th free module of (\ref{resolution}) is split off.

Now we are ready to establish the bound. Using  the inequalities (\ref{fact2}) and
(\ref{fact1}), we get for $i = 1,\dots,c+1$:
\begin{eqnarray*}
\lefteqn{M_1(J_1) \cdots M_i(J_1) \cdot [M_{i+1}(J_1) - d] \cdots [M_{c+1}(J_1) - d] } \\[1ex]
& = & M_1(J_1) \cdots M_{i-1} (J_1) \cdot [M_i(J_1) - d + d] \cdot [M_{i+1}(J_1) - d] \cdots
[M_{c+1}(J_1) - d]  \\[1ex]
& \geq & M_1(J_1) \cdots M_{i-1} (J_1) \cdot [M_{i}(J_1) - d] \cdots [M_{c+1}(J_1) - d] +
\\
& &
    d \cdot M_1(J_1) \cdots M_{i-1} (J_1) \cdot [M_{i+1}(J_1) - d] \cdots [M_{c+1}(J_1) - d] \\[1ex]
& \geq & M_1(J_1) \cdots M_{i-1} (J_1) \cdot [M_{i}(J_1) - d] \cdots [M_{c+1}(J_1) - d] +
\\
& &
    d \cdot M_1(I) \cdots M_{i-1} (I) \cdot M_{i} (I) \cdots M_{c}(I) \\[1ex]
& \geq & M_1(J_1) \cdots M_{i-1} (J_1) \cdot [M_{i}(J_1) - d] \cdots [M_{c+1}(J_1) - d] +
d \cdot c! \cdot e(R/I)
\end{eqnarray*}
where we also used the assumption on $I$. Applying this estimate repeatedly as well as
Condition (\ref{upper bd assumption}), we get
\begin{eqnarray*}
\lefteqn{ \frac{1}{(c+1)!} M_1(J_1) \cdots  M_{c+1}(J_1) } \\[1ex]
& \geq & \frac{1}{(c+1)!} M_1(J_1) \cdots M_{c} (J_1) \cdot [M_{c+1}(J_1) - d] +
\frac{d}{c+1}
\cdot e(R/I) \\[1ex]
& \geq & \frac{1}{(c+1)!} M_1(J_1) \cdots M_{c-1} (J_1) \cdot [M_c (J_1) - d] [
M_{c+1}(J_1) - d]  + \frac{2 d}{c+1}
\cdot e(R/I) \\[1ex]
& \geq & \ldots \\
& \geq & \frac{1}{(c+1)!} [M_1(J_1) - d] \cdots  [M_{c+1}(J_1) - d]  + d \cdot e(R/I) \\[1ex]
& \geq & \frac{1}{(c+1)!} M_1(J) \cdots  M_{c+1}(J) + d \cdot e(R/I) \\[1ex]
& \geq & e(R/J) + d \cdot e(R/I) \\[1ex]
& = & e(R/J_1),
\end{eqnarray*}
as claimed.
\end{proof}

\medskip

\begin{example}\label{all top deg split}
Let $R = k[x,y,z]$ and let
\[
\begin{array}{rcl}
 J & = & (x, y^9, z^6) \\
 I &  = & (y^9, z^6) \\
J_1 & = & (x^2,  y^9, z^6).
\end{array}
\]
Then $J_1$ is a basic double G-link as described above with $d=1$ and $L=x$, and we see
that
\[
\begin{array}{rclcrclcrcl}
M_1(J) & = & 9 && M_1(I) & = &  9 && M_1(J_1) & = & 9 \\
M_2(J) & = & 15 && M_2(I) & = & 15 && M_2(J_1) & = & 15 \\
M_3(J) & = & 16 && &&&& M_3(J_1) & = & 17.
\end{array}
\]
Furthermore, in all six cases there is only one summand of top degree in the
corresponding free  module.  Hence it follows that (\ref{upper bd assumption}) and
(\ref{upper bd assumption}$'$) are false in this example.
\end{example}

\medskip

\begin{remark}
We have not been able to find an example where Condition (\ref{fact2}) is not satisfied.
In fact, we suspect that this condition is always true.
\end{remark}

For later use, we record when we have equality in Theorem \ref{upper bd for BDL}. It is
an immediate consequence of its proof.

\begin{corollary} \label{cor-sharp-upper}
Adopt the notation and assumptions of Theorem \ref{upper bd for BDL}. Then the upper
bound of the Multiplicity conjecture is sharp for $J_1$ if and only if the upper bound of
the Multiplicity conjecture is sharp for $J$ as well as for $I$ and
$$
\begin{array}{rcll}
M_i (J_1) & = & M_i (J) + d, \quad  & i = 1,\ldots,c+1 \\
M_i (J_1) & =  & M_i (I), \quad & i = 1,\ldots,c \\
M_i (J_1) & =  & M_{i-1} (I) + d , \quad & i = 2,\ldots,c+1.
\end{array}
$$
\end{corollary}

For the lower bound of the Multiplicity Conjecture, there is a similar statement.

\begin{corollary} \label{lower bd for BDL}
Under the assumptions of Proposition \ref{KMMNP result}, assume that $I$ and $J$ both
 satisfy the lower bound of the Multiplicity conjecture \ref{conj1} (though we do not assume that
 $R/I$ nor $R/J$ are Cohen-Macaulay).
Assume further that
\begin{equation}\label{assm}
m_i(J_1) \leq m_{i-1}(I)+d \ \ \ \hbox{ for all $2 \leq i \leq c+1$}.
\end{equation}
Then $J_1$ also satisfies the lower bound of the Multiplicity conjecture.

Moreover, the lower bound of the Multiplicity conjecture is sharp for $J_1$ if and only
if the lower bound of the Multiplicity conjecture is sharp for $J$ as well as for $I$ and
$$
\begin{array}{rcll}
m_i (J_1) & = & m_i (J) + d, \quad  & i = 1,\ldots,c+1 \\
m_i (J_1) & =  & m_i (I), \quad & i = 1,\ldots,c \\
m_i (J_1) & =  & m_{i-1} (I) + d , \quad & i = 2,\ldots,c+1.
\end{array}
$$
\end{corollary}

\begin{proof}
The proof is almost identical to that of Theorem \ref{upper bd for BDL}.  This time we need information about the smallest shift in any free module.
An analysis similar to that in Theorem \ref{upper bd for BDL} shows immediately that
\begin{equation}\label{fact11}
m_i(J_1) \leq m_i(I), \ \ \  1 \leq i \leq c,
\end{equation}
since no summand of $\mathbb F_i$ in the $i$-th free module of
(\ref{resolution}) is split off.
  But furthermore, if a splitting of a summand $R(-t-d)$ between ${\mathbb F}_{i} (-d)$ (in the $(i+1)$-st free module) and ${\mathbb G}_i (-d)$ (in the $i$-th free module) occurs, then
$\mathbb F_i$ (in the $i$-th free module) contains the summand $R(-t)$, and $t
< t+d$.  Hence
\begin{equation} \label{fact22}
m_i(J_1) \leq m_i(J)+d, \ 1 \leq i \leq c.
\end{equation}
With the three inequalities (\ref{assm}), (\ref{fact11}) and (\ref{fact22}), the proof is almost identical to the proof of the upper bound in Theorem \ref{upper bd for BDL} and is left to the reader.
\end{proof}

The following consequence of Theorem \ref{upper bd for BDL} and Corollary \ref{lower bd
for BDL}  has been proven independently in \cite{H} using different methods.

\begin{corollary}
Let $J_1$ be the hypersurface section of $I$ by $F$, i.e.\ $J_1 = I + (F)$ where $F$ is a
homogeneous polynomial of degree $d > 0$ such that $I : F = I$. Let $c+1$ be the
codimension of $J_1$. Assume that the lower and upper bounds of the Multiplicity
conjecture hold for $I$. Then $J_1$ satisfies the conjectured bound, that is:
\[
\frac{1}{(c+1)!} m_1(J_1) m_2(J_1) \dots m_{c+1}(J_1) \leq e(R/J_1) \leq \frac{1}{(c+1)!}
M_1(J_1) M_2(J_1) \dots M_{c+1}(J_1)
\]
\end{corollary}

\begin{proof}
The proof is very similar to the proof above, but now there is no need for any extra
hypotheses.  Indeed, we begin with the exact sequence
\[
0 \rightarrow I(-d) \rightarrow I \oplus R(-d) \rightarrow J_1
\rightarrow 0.
\]
We again consider a minimal free resolution $\mathbb F_\bullet$ for $I$, and
now a mapping cone gives the following free resolution for $J_1$:
\[
0 \rightarrow  \mathbb F_c(-d) \rightarrow
\begin{array}{c}
\mathbb F_{c-1}(-d) \\
\oplus \\
\mathbb F_c
\end{array}
\rightarrow
\begin{array}{c}
\mathbb F_{c-2}(-d) \\
\oplus \\
\mathbb F_{c-1}
\end{array}
\rightarrow \cdots \rightarrow
\begin{array}{c}
\mathbb F_1(-d) \\
\oplus \\
\mathbb F_2
\end{array}
\rightarrow
\begin{array}{c}
\mathbb F_1 \\
\oplus \\
R(-d)
\end{array}
\rightarrow J_1\rightarrow 0
\]
This time, there is no possible splitting, so this is a minimal free resolution. For the
lower bound, we observe the (smaller) set of inequalities
\begin{equation} \label{facts}
\begin{array}{rcll}
m_i(J_1) & \leq & m_i(I) & \hbox{ for  } 1 \leq i \leq c \\
m_1(J_1) & \leq & d \\
m_i(J_1) & \leq & m_{i-1}(I) +d &  \hbox{ for } 2 \leq i \leq c+1
\end{array}
\end{equation}
Then almost the same proof as given above (again using the trick of adding and
subtracting $d$ several times) yields the desired bound.  It is
simpler since we
only have to ``convert" terms involving
$m_i(J_1)$ to terms involving $m_i(I)$; we do not have any $m_i(J)$ involved.
In the very last step we use the bound $m_1(J_1) \leq d$.  We omit the details.  The upper bound is proven similarly.
\end{proof}

Now we will discuss how Theorem \ref{upper bd for BDL} and Corollary  \ref{lower bd for
BDL} can be applied to show that the Multiplicity conjecture \ref{conj1} is true for
standard determinantal ideals.

Let $A$ be a homogeneous $t \times (t+c-1)$  matrix with entries in $R$, i.e.\ such that
multiplication by $A$ defines a graded homomorphism $\ffi: F \to G$ of free $R$-modules
and that all entries of degree zero are zero. Then we denote the ideal generated by the
maximal minors of $A$ by $I(\ffi) = I(A)$. Its codimension is at most $c$. We call $I
\subset R$ a {\em standard determinantal ideal} if $I = I(A)$ for some homogeneous $t
\times (t+c-1)$ matrix $A$ and $\codim I = c$.

\begin{corollary} \label{thm-det}
The Multiplicity conjecture  is true for all standard determinantal ideals.

Furthermore, the following conditions are equivalent for every standard determinantal
ideal $I = I(A)$ of codimension $c$:
\begin{itemize}
\item[(i)] $e(R/I) = \frac{1}{c!} \prod_{i=1}^c M_i$;
\item[(ii)] $e(R/I) = \frac{1}{c!} \prod_{i=1}^c m_i$;
\item[(iii)] $I$ has a pure minimal free resolution.
\item[(iv)] All the entries in $A$ have the same degree.
\end{itemize}
\end{corollary}

\begin{proof}
Our method of proof is essentially the same as the one in \cite{MNR} to show the claim in
codimension two.  This approach has also recently been carried out by Mir\'o-Roig
\cite{rosapaper} to establish the bounds, so we just give enough details to discuss
sharpness of the bounds. The result has also independently been shown by Herzog and Zheng
\cite{H}.

First, we discuss the bounds. We will switch to the notation of the basic double link
results. Let $J_1 = I(A) \subset R$ be a standard determinantal ideal defined by the
homogeneous $(t) \times (t+c-1)$ matrix $A$ ($t \geq 1)$ . We will show the claim by
induction on $t \geq 1$. If $t = 1$, then $J_1$ is a complete intersection and the bounds
are shown in \cite{HS}. Since, the bounds are trivial for principal ideals, we may assume
that $c \geq 2$ and that the claims are shown for ideals of codimension $\leq c-1$.

Let $t \geq 2$. By reordering rows and columns we can arrange that the degree of the
entries of  $A$ increase from bottom to top and from left to right. Let us write the
resulting degree matrix of $A$ as follows:
\begin{eqnarray} \label{deg mat}
\partial A & = & \begin{bmatrix}
a_{1 1} & a_{1 2} & \ldots & a_{1 c} & & *\\
& a_{2 1} & a_{2 2} & \ldots & a_{2 c} & & \\
& & \ddots & \ddots &  &  \ddots\\
* & & & a_{t 1} & a_{t 2} & \ldots & a_{t c} \\
\end{bmatrix}
\end{eqnarray}
Note that our ordering of degrees means that
\begin{equation} \label{order}
a_{i-1, j}\leq a_{i j} \leq a_{i, j+1}.
\end{equation}
Observe that the whole degree matrix $\partial A$ is completely determined by the entries
specified in (\ref{deg mat}). Moreover, the Eagon-Northcott complex shows that the graded
Betti numbers of $J_1$ are completely determined by $\partial A$. Hence, it suffices to
show the bounds for $J_1 = I(A)$ where
\begin{eqnarray} \label{mon-mat}
A & = & \begin{bmatrix}
x_1^{a_{1 1}} & x_2^{a_{1 2}} & \ldots & x_c^{a_{1 c}} & & 0\\
& x_1^{a_{2 1}} & x_2^{a_{2 2}} & \ldots & x_c^{a_{2 c}} & & \\
& & \ddots & \ddots &  &  \ddots\\
0 & & & x_1^{a_{t 1}} & x_2^{a_{t 2}} & \ldots & x_c^{a_{t c}} \\
\end{bmatrix}
\end{eqnarray}
Denote by $B$ the matrix that is obtained from $A$ by deleting the last column and by
$A'$ the matrix that one gets after deleting the last row of $B$. Then $J := I(A')$ and
$I := I(B)$ are standard determinantal ideals satisfying the relation
$$
J_1 = x_c^{a_{t c}} \cdot J + I,
$$
where $I : x_c^{a_{t c}} = I$, i.e.\ $J_1$ is a basic double link of $J$.

The Eagon-Northcott complex provides that
\begin{eqnarray*}
m_1 (J_1) & = & a_{1 1} + a_{2 1} + \ldots + a_{t 1} \\
m_i (J_1) & = & m_{i-1} (J_1) + a_{t i}, \quad i = 2,\ldots,c \\
M_1(J_1) & = & a_{1 c} + a_{2 c} + \ldots + a_{t c} \\
M_i (J_1) & = & M_{i-1} (J_1) + a_{t, c+1-i}, \quad i = 2,\ldots,c
\end{eqnarray*}
and similar formulas for $J$ and $I$. Hence, using our ordering (\ref{order}), it is easy
to see that the assumptions of Theorem \ref{upper bd for BDL} and Corollary  \ref{lower
bd for BDL} are satisfied, thus it follows that the bounds are true for $J_1$.

Now we discuss sharpness. It suffices to show that either condition (i) or (ii) implies
(iv). Assume (i) is true. Then, Corollary \ref{cor-sharp-upper} and the induction
hypothesis show that all the entries of $B$ have the same degree. But Corollary
\ref{cor-sharp-upper} also provides $M_{i-1} (I) = M_i (J_1)$, which forces the entries
in the last column of $A$ to have the same degree as the entries of $B$.

The proof that (ii) implies (iv) is similar; we omit the details.
\end{proof}

As a further application of basic double links, we conclude this section with a closed
formula for the degree of a standard determinantal ideal.

\begin{theorem} \label{thm-degree-form}
Let $\ffi: \oplus_{j = 1}^{t+c-1} R(-d_j) \to F := \oplus_{i = 1}^t R(-b_i)$ be a
homomorphism of graded $R$-modules such that $\im \ffi \subset \fm \cdot F$ and $\coker
\ffi$ has codimension $c \geq 1$. Assume that $b_1 \leq b_2 \leq \ldots b_t$ and $d_1
\leq d_2 \leq \ldots \leq d_{t+c-1}$. Then the ideal of maximal minors of $\ffi$ has
degree
$$
\deg I(\ffi) = \sum_{i_c = 1}^t a_{i_c, c} \left \{ \sum_{i_{c-1} = 1}^{i_c} a_{i_{c-1},
c-1} \left \{ \ldots \left \{ \sum_{i_2 = 1}^{i_3} a_{i_2, 2} \cdot \prod_{i_1 = 1}^{i_2}
a_{i_1, 1} \right \} \ldots \right \} \right \}
$$
where $a_{i j} := d_{i+j-1} - b_i$, $i = 1,\ldots,t$ and $j= 1,\ldots,c$ and the formula
reads for $c=1$ as $\deg I(\ffi) = \prod_{i_1 = 1}^{t} a_{i_1, 1}$.
\end{theorem}

\begin{proof}
We use the approach employed in the proof of Corollary \ref{thm-det}. We may assume that
the map $\ffi$ is represented by the matrix $A$ as specified in (\ref{mon-mat}). Denoting
the matrices $B$ and $A'$ obtained from $A$ by deleting the last column and row as in the
proof of Corollary \ref{thm-det}, we get
$$
I(A) = x_c^{a_{t c}} \cdot I(A') + I(B),
$$
thus
$$
\deg I(A) = a_{t c} \cdot  \deg I(A') + \deg I(B).
$$
The cases $c=1$ and $t=1$ being trivial, our claim follows now easily by induction on $c$
and $t$.
\end{proof}

\begin{remark}
In the classical case $b_1 = \ldots = b_t = 0$ and $d_1 = \ldots = d_{t+c-1} = 1$, i.e.\
all entries of $A$ have degree one,  Theorem \ref{thm-degree-form} specializes to
$$
\deg I(\ffi) = \binom{t+c-1}{c},
$$
which of course also follows by  the classical Porteous' formula.
\end{remark}

\section{Modules with quasi-pure resolutions} \label{sec-quasi-pure}
Herzog and Srinivasan  proved in \cite{HS} Section 1 the Multiplicity conjecture for
$k$-algebras with quasi-pure resolutions. In this section we generalize this result to
the module case.  Let $N$ be a finitely generated graded $R$-module. Let
$$
F_{\pnt}: 0 \to \bigoplus_{j=1}^{b_p} R(-d_{pj}) \to \dots \to \bigoplus_{j=1}^{b_0}
R(-d_{0j}) \to 0
$$
be the minimal graded free resolution of $N$ and define the invariants $M_i$ and $m_i$ as
in the Introduction. We say that $N$ has a {\em quasi-pure resolution} if $m_i\geq
M_{i-1}$ for all $1\leq i \leq p$. We follow the line of proof in \cite{HS} with the
necessary changes.

\begin{lemma}
\label{quasihelper} Let  $N$ be a graded $R$-module of codimension  $c\geq 1$. Then
$$
\sum_{i=0}^p (-1)^i \sum_{j=1}^{b_i} d_{ij}^k
=
\begin{cases}
0               & \text{if } 1\leq k < c,\\
(-1)^c c! \cdot e(N)  & \text{if } k=c.
\end{cases}
$$
\end{lemma}
\begin{proof}
Let $H_N(t)$ be the Hilbert series of $N$. Since $N$ has
a quasi-pure resolution, we can compute this series as
$$
H_N(t)
=
\frac{\sum_{i=0}^p \sum_{j=1}^{b_i} (-1)^it^{d_{ij}} }{(1-t)^n}.
$$
On the other hand we know that
$$
H_N(t) = \frac{Q(t)}{(1-t)^d}
$$
where $Q(t)$ is a polynomial such that $Q(1)=e(N)$ equals the multiplicity of $N$
and $d=\dim N$.
Thus we get that
$$
\sum_{i=0}^p \sum_{j=1}^{b_i} (-1)^i t^{d_{ij}}
=
Q(t) (1-t)^c.
$$
Let $P(t)$ be the polynomial on the left hand side of this equation.
It follows that
\begin{eqnarray}
\label{nicetohave}
P^{(k)}(1)
=
\begin{cases}
0 & \text{for } 1\leq k < c,\\
(-1)^c c! \cdot e(N)& \text{for } k=c.
\end{cases}
\end{eqnarray}
We prove by induction
on $k\in \{1,\dots, c\}$ and another induction on $l \in \{1,\dots, k\}$
that
$$
P^{(k)}(1)
=
\sum_{i=0}^p \sum_{j=1}^{b_i} (-1)^i d_{ij}^l(d_{ij}-1) \dots (d_{ij}-k+1 +l-1 ).
$$
The cases for $l=k$ together with (\ref{nicetohave}) give the desired formula of the lemma.
For $k=l=1$ and $k>1$, $l=1$ we have by the definition of $P(t)$ that
$$
P^{(k)}(1)
=
\sum_{i=0}^p \sum_{j=1}^{b_i} (-1)^i d_{ij}(d_{ij}-1) \dots (d_{ij}-k+1).
$$
Observe that we do not have to distinguish whether the $d_{ij}$
are bigger or equal to $k$ since we can add the remaining terms which are zero.
Assume that $k>1$ and $l>1$ and we know by the second induction hypothesis for $l$
that
$$
P^{(k)}(1)
=
\sum_{i=0}^p \sum_{j=1}^{b_i} (-1)^i d_{ij}^{l-1}(d_{ij}-1) \dots (d_{ij}-k+1 +(l-1)-1 ).
$$
It follows from the induction hypothesis on $k$ that
$$
P^{(k-1)}(1)
=
\sum_{i=0}^p \sum_{j=1}^{b_i} (-1)^i d_{ij}^{l-1}(d_{ij}-1) \dots (d_{ij}-k + 2 + (l-1)-1 ).
$$
On the other hand we know $P^{(k-1)}(1)=0$ as was already observed in (\ref{nicetohave}).
Hence we have that
$$
P^{(k)}(1)
=
P^{(k)}(1)
+
(k -l+1)
P^{(k-1)}(1)
=
\sum_{i=0}^p \sum_{j=1}^{b_i} (-1)^i d_{ij}^l(d_{ij}-1) \dots (d_{ij}-k+1 +l-1 )
$$
as desired.
\end{proof}

For modules with quasi-pure resolution we get the following result, which hints at a
possible generalization of the Multiplicity conjecture to the case of modules.

Notice the degree of a maximal module is simply equal to its rank. Thus, we exclude this
case and focus on torsion modules throughout the remainder of this note.

\begin{theorem}
\label{modulequasi} Let $N=\dirsum_{i\in\Z}N_i$ be a finitely generated graded $R$-module
with $p=\projdim N$. Assume that $N$ is Cohen-Macaulay, $\rank N=0$, and that $N$ has a
quasi-pure resolution. Then
$$
\prod_{i=1}^c (m_i-M_0) / p! \leq e(N) \leq \prod_{i=1}^c (M_i - m_0)/p!
$$
with equality below (resp.\ above) if and only if $N$ has a pure resolution.
\end{theorem}
\begin{proof}
We consider the $p+1\times p+1$-square matrix
$$
A
=
\left(
  \begin{array}{ccccc}
    b_0                       & b_1                       & \cdots      & b_p \\
    \sum_{j=1}^{b_0} d_{0j}   & \sum_{j=1}^{b_1} d_{1j}   & \cdots & \sum_{j=1}^{b_p} d_{pj}   \\
    \vdots                    & \vdots                    & \vdots & \vdots \\
    \sum_{j=1}^{b_0} d_{0j}^p & \sum_{j=1}^{b_1} d_{1j}^p & \cdots & \sum_{j=1}^{b_p} d_{pj}^p   \\
  \end{array}
\right)
$$
By replacing the first column of $A$ by the alternating sum of all columns of $A$, we
obtain a matrix $A'$ such that $ \det (A) = \det(A')$. Since $\sum_{i=0}^p b_i=\rank N=0$
and, by Lemma \ref{quasihelper},  $ \sum_{i=0}^p (-1)^i \sum_{j=1}^{b_i} d_{ij}^p=(-1)^p
p! e(N), $ it follows that the first column of $A'$ is the transpose of the vector
$(0,\dots,0,(-1)^p p! e(N))$. Hence by expanding the determinant of $A'$ with respect to
the first column, we get
$$
\det(A)=\det(A')= p! e(N) \det(B),
$$
where $B$ is the $p\times p$-matrix
$$
B
=
\left(
  \begin{array}{ccccc}
    b_1                       & b_2                       & \cdots      & b_{p} \\
    \sum_{j=1}^{b_1} d_{1j}   & \sum_{j=1}^{b_2} d_{2j}   & \cdots & \sum_{j=1}^{b_{p}} d_{pj}   \\
    \vdots                    & \vdots                    & \vdots & \vdots \\
    \sum_{j=1}^{b_1} d_{1j}^{p-1} & \sum_{j=1}^{b_2} d_{2j}^{p-1} & \cdots & \sum_{j=1}^{b_{p}} d_{pj}^{p-1}   \\
  \end{array}
\right)
.
$$
Hence
\begin{eqnarray}
\label{eqeins}
\det(A)= p! e(N) \det(B)
=
p! e(N)
\sum_{j_1=1}^{b_1}
\sum_{j_2=1}^{b_2}
\dots
\sum_{j_p=1}^{b_p}
U(d_{1j_1}, \dots, d_{p j_p})
\end{eqnarray}
with the Vandermonde determinants
$$
U(d_{1j_1},\dots,d_{pj_p})
=
\det
\left(
  \begin{array}{ccccc}
    1          & 1          & \cdots   & 1 \\
    d_{1j_1}   & d_{2j_2}   & \cdots   & d_{pj_p}   \\
    \vdots                    & \vdots & \vdots & \vdots \\
    d_{1j_1}^{p-1} & d_{2j_2}^{p-1} & \cdots & d_{pj_p}^{p-1}   \\
  \end{array}
\right)
.
$$
We may also directly compute the determinant of $A$ as
\begin{eqnarray}
\label{eqzwei}
\det(A)
=
\sum_{j_0=1}^{b_0}
\sum_{j_1=1}^{b_1}
\dots
\sum_{j_p=1}^{b_p}
V(d_{0j_0},\dots,d_{pj_p})
\end{eqnarray}
with the corresponding Vandermonde determinants
$$
V(d_{0j_0},\dots,d_{pj_p})
=
\det
\left(
  \begin{array}{ccccc}
    1          & 1          & \cdots   & 1 \\
    d_{0j_0}   & d_{1j_1}   & \cdots   & d_{pj_p}   \\
    \vdots                    & \vdots & \vdots & \vdots \\
    d_{0j_0}^p & d_{1j_1}^p & \cdots & d_{pj_p}^p   \\
  \end{array}
\right)
.
$$
Since $N$ has a quasi-pure resolution we have that $d_{ij} \geq d_{i-1 k}$ for all $i,j,k$
and thus all the involved Vandermonde determinants are always non-negative.
Observe that
$$
V(d_{0j_0},\dots,d_{pj_p})
=
\prod_{i=1}^p (d_{ij_i} - d_{0j_0})
U(d_{1j_1},\dots,d_{pj_p})
$$
and thus
$$
\prod_{i=1}^p (m_i-M_0)\cdot
U(d_{1j_1},\dots,d_{pj_p})
\leq
V(d_{0j_0},\dots,d_{pj_p})
\leq
\prod_{i=1}^p (M_i-m_0)
\cdot
U(d_{1j_1},\dots,d_{pj_p})
$$
with equality below or above for all involved indices
if and only if $N$ has a pure resolution.
Hence it follows from (\ref{eqeins}) and (\ref{eqzwei})
that
$$
\prod_{i=1}^p (m_i-M_0)
\leq
p!e(N)
\leq
\prod_{i=1}^p (M_i-m_0)
$$
with equality below or above
if and only if $N$ has a pure resolution.
This concludes the proof.
\end{proof}

For modules with pure resolution we get the following nice formula for the multiplicity.
It generalizes the Huneke-Miller formula for the multiplicity of ideals with a pure
resolution \cite{HM}.

\begin{corollary}
Let $N=\dirsum_{i\in \N}N_i$ be a finitely generated graded $R$-module of codimension
$c$. Assume that $N$ is Cohen-Macaulay, $\rank N=0$, and that $N$ has a pure resolution $
 0 \to R(-d_{c})
\to \dots \to R(-d_{0}) \to 0$. Then
$$
e(N)=\prod_{i=1}^c (d_i-d_0) / c!
$$
\end{corollary}

\section{Modules of codimension 2} \label{sec-mod-2}

In \cite{roemer} the third author proved Conjecture \ref{conj2} for $k$-algebras of
codimension 2. The goal of this section will be to generalize this result to the module
case. As in \cite{roemer}, one of our tools is the use of general hyperplane sections.
Let $N=\dirsum_{i\in \N}$ be a finitely generated graded $R$-module. Following
\cite{ARHE00} we call an element $x \in R_{1}$ {\em almost regular} for $N$ if
$$
(0:_N x)_{a} = 0 \text{ for } a \gg 0.
$$
A sequence $x_1,\ldots,x_t \in R_{1}$ is an {\em almost regular sequence} for $N$ if for
all $i \in \{1,\dots,t\}$ the element $x_i$ is almost regular for
$N/(x_1,\ldots,x_{i-1})N$. Since neither the Betti numbers nor the multiplicity of $N$
changes by enlarging the field, we  may always assume that $k$ is an infinite field. It
is well-known that then after a generic choice of coordinates we can achieve that a
$k$-basis of $R_{1}$ is almost regular for $N$.

Let $x \in R_1$ be generic and almost regular for $N$. We observe that the following
holds:
\begin{enumerate}
\label{almosthelper}
\item
If $\dim N>0$, then $\dim N/xN=\dim N -1$.
\item
If $\dim N>1$, then $e(N) = e(N/xN)$.
\item
If $\dim N = 1$, then $e(N)\leq e(N/xN)$
with equality if and only if $N$ is Cohen-Macaulay.
\end{enumerate}
We are ready to prove the main theorem of this section.

\begin{theorem}
\label{modulecodim2} Let $N=\dirsum_{i\in \Z}$ be a finitely generated graded $R$-module
with $\rank N=0$.
\begin{enumerate}
\item
If $\codim N =1$, then
$$
e(N) \leq M_1-m_0
$$
with equality if and only if $N$ is Cohen-Macaulay and $N$ has a pure resolution.
\item
If $\codim N =2$, then
$$
e(N) \leq (M_1-m_0) (M_2-m_0)/2
$$
with equality if and only if $N$ is Cohen-Macaulay and $N$ has a pure resolution.
\end{enumerate}
\end{theorem}
\begin{proof}
We only prove (b) since the proof of (a) is simpler and is shown analogously. As noticed
above we may assume that  $x_1,\ldots,x_n \in R_1$ is a generic almost regular sequence
for $N$.

Let $\xb=x_1,\ldots,x_{n-2}$ and consider $\Tilde{N}=N/\xb N$. Observe that $2=\codim
N=\codim \Tilde{N}$ and $e(N)\leq e(\Tilde{N})$ with equality if and only if $N$ is
Cohen-Macaulay. Notice that $\Tilde{N}$ is a finitely generated graded
$\Tilde{R}$-module, where $\Tilde{R}$ is the 2-dimensional polynomial ring $R/\xb R$. Let
$$
M_i=\max\{j\in \Z\colon \beta^R_{i,j}(N)\neq 0\} \text{ for } i=1,2
$$
and
$$
\Tilde{M}_i=\max\{j\in \Z\colon \beta^{\Tilde{R}}_{i,j}(\Tilde{N})\neq 0\} \text{ for }
i=1,2.
$$
We claim that
\begin{equation}
\label{claim}
\Tilde{m_0}=m_0,\ \Tilde{M}_1 \leq M_1 \text{ and } \Tilde{M}_2\leq M_2.
\end{equation}
Note that for a Cohen-Macaulay module $N$ we have equalities everywhere.

Since $\dim\Tilde{N} =0$, the module $\Tilde{N}$ is Cohen-Macaulay. Hence
$\Tilde{M}_2>\Tilde{M}_1>\Tilde{M}_0$ and $\Tilde{m}_2>\Tilde{m}_1>\Tilde{m}_0$. Let
$\Tilde{N} \cong \Tilde{F}/\Tilde{G}$ where $\Tilde{F}$ is the first finitely generated
graded free $\Tilde{R}$-module of a minimal free resolution of $\Tilde{N}$ and
$\Tilde{G}$ the kernel of the map $\Tilde{F} \to \Tilde{N}$. Consider the module
$\Tilde{N}'=\Tilde{F}/\Tilde{G}_{\geq M_1}$. Observer that $\Tilde{N}'$ is still
Cohen-Macaulay of codimension 2 because $\dim \Tilde{N}'=0$. Let
$\Tilde{M}_i',\Tilde{m}_i'$ denote the corresponding invariants of $\Tilde{N}'$. We see
that
$$
\Tilde{M}_0'=\Tilde{M}_0, \
\Tilde{m}_0'=\Tilde{m}_0,\
\Tilde{M}_1'=\Tilde{m}_1'=\Tilde{M}_1,
\text{ and }
\Tilde{M}_2'=\Tilde{M}_2.
$$
Furthermore we have that
$$
\Tilde{m}_2'> \Tilde{m}_1'=\Tilde{M}_1'.
$$
Thus $\Tilde{N}'$ is Cohen-Macaulay of rank $0$ and has a quasi-pure resolution. Hence we
may apply Theorem \ref{modulequasi} to the module $\Tilde{N'}$. Note that
$$
e(\Tilde{N}) \leq e(\Tilde{N}')
$$
with equality if and only if $\Tilde{M}_1=\Tilde{m}_1$.
All in all we get that
\begin{eqnarray*}
e(N)
&\leq& e(\Tilde{N})\\
&\leq& e(\Tilde{N}')\\
&\leq& (\Tilde{M}_1'-\Tilde{m}_0')(\Tilde{M}_2'-\Tilde{m}_0')\\
& = & (\Tilde{M}_1-\Tilde{m}_0)(\Tilde{M}_2-\Tilde{m}_0)\\
& \leq & (M_1-m_0)(M_2-m_0)
\end{eqnarray*}
with equalities everywhere if and only if  $N$ is Cohen-Macaulay and has a pure resolution.

It remains to prove claim (\ref{claim}). The first two inequalities can easily be seen:
$\Tilde{m}_0$ is the minimal degree of a minimal generator of $\Tilde{N}$ and $m_0$ is
the minimal degree of a minimal generator of $N$. By Nakayama's lemma this degree does
not change by passing from $N$ to $\Tilde{N}$. If $N=F/G$ where $F$ is the first finitely
generated graded free $R$-module of a minimal free resolution of $N$ and $G$ the kernel
of the map $F \to N$, then $\Tilde{N}\cong \bigl(F/\xb F\bigr)/ \bigl(G+\xb F/\xb
F\bigr)$. Thus we see that $\Tilde{M}_1 \leq M_1$.

To prove the remaining inequality $\Tilde{M}_2\leq M_2$, we can use the proof of Theorem
2.4 in \cite{roemer} word by word since that proof holds also in the module case. This
completes the argument.
\end{proof}

\begin{remark}
Theorem \ref{modulecodim2} extends Theorem 2.4 in \cite{roemer} (cf.\ also \cite{H},
Theorem 3.1) from cyclic modules to arbitrary modules (of codimension two). It proves
Conjecture \ref{conj2} for modules whose codimension is at most two.
\end{remark}

Theorem \ref{modulequasi} suggests for a Cohen-Macaulay torsion module $N$ of codimension
$c$ that $\frac{1}{c!} \prod_{i=1}^c (m_i - M_0) \leq e(N)$. While this would be an
interesting bound in some cases, for example, if all generators of $N$ have the same
degree, it does not always give useful information because this lower bound can be a
negative number.



\begin{thebibliography}{999}
\bibitem{ARHE00}
A.\ Aramova and J.\ Herzog, {\it Almost regular sequences and Betti numbers}. Amer.\ J.\
Math.\ {\bf 122} (2000), 689--719.

\bibitem{macaulay} D.\ Bayer and M.\ Stillman, Macaulay: A system for
computation in algebraic geometry and commutative algebra. Source and object code
available for Unix and Macintosh computers.  Contact the authors, or download from
ftp://math.harvard.edu via anonymous ftp.

\bibitem{F}
C.\ Francisco, {\em New approaches to bounding the multiplicity of an ideal}, Preprint,
2005.

\bibitem{FS} C.\ Francisco and H.\ Srinivasan, {\em Multiplicity conjectures}, Preprint, 2005.

\bibitem{gold} L.H.\ Gold, {\em A degree bound for codimension two lattice
ideals}, J.\ Pure Appl.\ Algebra {\bf 183} (2003), 201--207.

\bibitem{GSS} L.H.\ Gold, H.\ Schenck and H.\ Srinivasan, {\em Betti
numbers and degree bounds for some linked zero-schemes}, Preprint, 2004.

\bibitem{GV} E.\ Guardo and A.\ Van Tuyl, {\em Powers of complete
intersections: graded Betti numbers and applications}, to appear in Ill.\ J.\ Math.

\bibitem{HS}
J.\ Herzog and H.\ Srinivasan, {\em Bounds for multiplicities}. Trans.\ Amer.\ Math.\
Soc. {\bf 350} (1998), no. 7, 2879--2902.

\bibitem{H}
J.\ Herzog  and Zheng, {\em Notes on the multiplicity conjecture}, Preprint, 2005.

\bibitem{HM} C.\ Huneke and M.\ Miller, {\em A note on the multiplicity of
Cohen-Macaulay algebras with pure resolutions}, Canad.\ J.\ Math.\ {\bf 37} (1985),
1149--1162.

\bibitem{HU} C.\ Huneke and B.\ Ulrich, {\em General Hyperplane Sections of
Algebraic Varieties}, J.\ Alg.\ Geom.\ {\bf 2} (1993), 487--505.

\bibitem{KMMNP} J.\ Kleppe, J.\ Migliore, R.M.\ Mir\'o-Roig, U.\ Nagel, and C.\
Peterson, {\em Gorenstein Liaison, Complete Intersection Liaison Invariants and
Unobstructedness}, Memoirs of the Amer.\ Math.\ Soc.\ {\bf 154}, 2001; 116 pp.

\bibitem{migbook} J.\ Migliore, ``Introduction to Liaison Theory and
Deficiency Modules,''   Progress in Mathematics {\bf 165}, Birkh\"auser, 1998.

\bibitem{london} J.\ Migliore, {\em Submodules of the deficiency module}, J.\
Lond.\ Math.\  Soc.\ {\bf  48}(3) (1993), 396--414.

\bibitem{MN4} J.\ Migliore and U.\ Nagel, {\em Monomial Ideals and the
Gorenstein Liaison Class of a Complete Intersection}, Compositio Math.\ {\bf
133} (2002), 25--36.

\bibitem{MNR} J.\ Migliore, U.\ Nagel and T.\ R\"omer, {\em The Multiplicity
Conjecture in low codimensions}, to appear in Math.\ Res.\ Lett.

\bibitem{rosapaper} R.\ Mir\'o-Roig, {\em A note on the multiplicity of determinantal ideals},
Preprint, 2005.

\bibitem{rao} P.\ Rao, {\em Liaison among Curves in} $\PP^3$, Invent.\ Math.\
{\bf 50} (1979), 205--217.

\bibitem{roemer} T.\ R\"omer, {\em Note on bounds for multiplicities}, J.\ Pure Appl.\
Algebra  {\bf 195}  (2005), 113--123.

\bibitem{srinivasan} H.\ Srinivasan, {\em A note on the multiplicities of
Gorenstein algebras}, J.\ Algebra {\bf 208} (1998), no. 2, 425--443.


\end{thebibliography}
\end{document}